\documentclass[12pt]{amsart}%
\usepackage{amsfonts}
\usepackage{amsmath}
\usepackage{amsthm, amscd}
\usepackage[v2]{xy}
\usepackage{amssymb}
\usepackage{graphicx}
\usepackage{enumerate}
\usepackage{color}
\usepackage{setspace}

\newcommand{\dd}{\,\textrm{d}}

\newtheorem{theorem}{Theorem}[section]
\newtheorem{corollary}[theorem]{Corollary}
\newtheorem{lemma}[theorem]{Lemma}
\newtheorem{proposition}[theorem]{Proposition} 
\theoremstyle{definition}
\newtheorem{definition}[theorem]{Definition}

\numberwithin{equation}{section}
\theoremstyle{remark}  
\newtheorem{remark}[theorem]{Remark}

\def\R{{\mathbb R}}

\def\Z{{\mathbb Z}}
\def\N{{\mathbb N}}

\begin{document}
\title{An operator-fractal}
\author[P.E.T. Jorgensen]{Palle E. T. Jorgensen}
\address[Palle E.T. Jorgensen]{Department of Mathematics, The University of Iowa, Iowa
City, IA 52242-1419, U.S.A.}
\email{jorgen@math.uiowa.edu}
\urladdr{http://www.math.uiowa.edu/\symbol{126}jorgen/}
\author[K.A. Kornelson]{Keri A. Kornelson}
\address[Keri Kornelson]{Department of Mathematics, The University of Oklahoma, Norman, OK, 73019-0315, U.S.A.}
\email{kkornelson@ou.edu}
\urladdr{http://www.math.ou.edu/\symbol{126}kkornelson/}
\author[K.L. Shuman]{Karen L. Shuman}
\address[Karen Shuman]{Department of Mathematics and Statistics,
Grinnell College, Grinnell, IA 50112-1690, U.S.A.}
\email{shumank@grinnell.edu}
\urladdr{http://www.math.grinnell.edu/\symbol{126}shumank/}

\thanks{The second and third authors were supported in part by NSF grant DMS-0701164.}



\begin{abstract}
Certain Bernoulli convolution measures $\mu$ are known to be spectral.  Recently, much work has concentrated on determining conditions under which orthonormal Fourier bases (i.e. spectral bases) exist.  For a fixed measure known to be spectral, the ONB need not be unique; indeed, there are often families of such spectral bases. 

Let $\lambda = \frac{1}{2n}$ for a natural number $n$ and consider the Bernoulli measure with scale factor $\lambda$.  It is known that $L^2(\mu_{\lambda})$ has a Fourier basis.  We first show that there are Cuntz operators acting on this Hilbert space which create an orthogonal decomposition, thereby offering powerful algorithms for computations for Fourier expansions.

When $L^2(\mu_{\lambda})$ has more than one Fourier basis, there are natural unitary operators $U$, indexed by a subset of odd scaling factors $p$; each $U$ is defined by mapping one ONB to another.   We show that the unitary operator $U$ can also be orthogonally decomposed according to the Cuntz relations.  Moreover, this \textit{operator-fractal} $U$ exhibits its own self-similarity. 
\end{abstract}
\maketitle

\tableofcontents

\thispagestyle{empty}

\section{Introduction}

In the theory of fractals and fractal measures, patterns of symmetry and self-similarity abound.  In this paper, we study a new class of unitary operators on fractal $L^2$ spaces.  This class of operators have surprising self-similar properties which, to our knowledge, have not been observed before.  In one special case, we discover even more intricate patterns in the matrix of the operator. 

Fractal sets which are invariant under a collection of contractive maps (an iterated function system) exhibit scaling self-similarity; one sees the same shape when zooming in to look more closely at one part of the set.   We might call this scaling \textit{in the small} to find self-similarity.  It is also common to find that the discrete sets which index Fourier bases for fractal $L^2$ spaces have an expansive self-similarity themselves.  We say that these sets (called \textit{spectra}) have a self-similarity \textit{in the large}.   When a measure $\mu$ has such a discrete index set $\Gamma$ for a Fourier basis, we call $(\mu, \Gamma)$ a \textit{spectral pair}.
 
For illustration, consider first the simplest case---the Bernoulli convolution formed by recursive scaling by $\frac14$ with two affine maps on the real line. The resulting measure $\mu_{\frac14}$ is an infinite Bernoulli convolution, also called a Cantor measure, or a Hutchinson measure.  The Hilbert space $L^2(\mu_{\frac14})$ has a Fourier basis which has a self-similarity under scaling in the large by $4$.    It is somewhat surprising to find that that $\mu_{\frac14}$ also gives rise to a symmetry based on scaling by $5$.  It turns out that scaling by $5$ transforms the Fourier spectrum $\Gamma$ into another ONB $5\Gamma$. As a result, we have a natural unitary operator $U$ acting in $L^2(\mu_{\frac14})$ which maps one ONB to the other.  Its spectral properties turn out to reveal a surprising level of symmetry and self-similarity which lead us to the nomenclature \textit{operator-fractal}.  
  
After briefly reviewing Bernoulli convolutions and some of their properties in Section \ref{sec:background}, we demonstrate the intimate connection between the spectral theory of affine fractals and representations of the Cuntz algebras from the theory of $C^*$-algebras in Section \ref{sec:Cuntz}.  In particular, the scaling operations on the spectrum $\Gamma$ influence the structure of the Hilbert space $L^2(\mu)$.  We find that the scaling operators are isometries which satisfy Cuntz relations (Definition \ref{Defn:CuntzRel}); the Cuntz isometries induce the orthogonal decomposition of $L^2(\mu)$ shown in Proposition \ref{Prop:decomp}.   

For certain odd integers $p$, the scaled set $p\Gamma$ also is a spectrum for a Fourier basis on $L^2(\mu)$ (see \cite{LaWa02, DHS09, DJ09, JKS11a} for details).  We denote by $U$ the unitary operator which enacts the scaling by $p$ on the ONB $\Gamma$.  In Section \ref{sec:U}, we examine the interplay between $U$ and the isometries that form the Cuntz algebra.  In Section \ref{Subsec:Commute}, we prove commutation relations for $U$ and one of the Cuntz operators $S_0$.  More applications of Cuntz algebras to the setting of fractal measures can be found in \cite{BrJo99}.
 
In our main theorem, Theorem \ref{thm:block}, we show that when $(\frac{1}{2n}, p\Gamma)$ is a spectral pair, then the associated operator $U$ is an orthogonal sum of a rank-one projection and infinitely many copies of the same operator; these copies arise from the Cuntz decomposition.  In the special $(4,5)$ case, the operator which appears infinitely many times in $U$ is just a unitary multiplication operator on $U$ itself (Theorem \ref{Thm:Multiplication}).  Finally, we study the structure of the matrix of $U$ acting on the first Cuntz subspace in the $(4,5)$ case (Theorem \ref{Thm:SpecialCase}).    
  
  \section{Background and Notation}\label{sec:background}

In this section, we outline the construction of the Bernoulli convolution measures and provide some of the necessary background results.   

\subsection{Affine iterated function systems}\label{Subsec:IFS}

The \textit{Bernoulli convolution measures} are a special class of fractal measures, those generated on the real line by two affine maps with scale factor $\lambda \in (0,1)$.  There is a rich literature on these measures; see \cite{Erd39, PeSo96, PSS00, Sid03, LaWa02, JKS07b, JKS08, DJ09, DHS09, DHSW11, JKS11a} for just a small sampling.  The problem we consider here fits into a wider framework of spectral pairs, self-similarity, and spectral duality.  For the interested reader, other related work in the literature includes \cite{Fug74, Tao04, GaNa98, GaYu06, PeWa01, Wan02, Li07}.  The interplay between the scale factor $\lambda$ and a spectral scaling between orthonormal bases has been considered in \cite{JR95}.

A Bernoulli convolution measure $\mu_{\lambda}$ on the real line can be realized as an invariant measure for an iterated function system (IFS).  Let $\{\tau_+, \tau_-\}$ be the IFS \begin{equation} \tau_+(x) = \lambda(x+1), \qquad \tau_-(x) = \lambda(x-1),\end{equation} where the scaling factor $\lambda$ is in  $(0,1)$.  

By the well-established algorithm from \cite{Hut81}, the IFS $\{\tau_+, \tau_-\}$ generates a compact \textit{attractor set}, $X_{\lambda}$,  that satisfies 
\[ X_{\lambda} = \tau_+(X_{\lambda}) \cup \tau_-(X_{\lambda}).\]  
The IFS also generates a measure $\mu_{\lambda}$ (the Bernoulli convolution), supported on $X_{\lambda}$, which satisfies a similar invariance equation: 
\begin{equation}\label{eqn:mu_invariance} 
\mu_{\lambda} = \frac12 \Bigl(\mu_{\lambda} \circ \tau_+^{-1}\Bigr) + \frac12 \Bigl(\mu_{\lambda} \circ \tau_-^{-1}\Bigr).
\end{equation}   

We might ask under which conditions the measure $\mu_{\lambda}$ has a Fourier basis (i.e. an orthonormal basis of complex exponential functions) for the Hilbert space $L^2(\mu_{\lambda})$.  When such an orthonormal Fourier basis exists, we say that $\mu_{\lambda}$ is a \textit{spectral measure}.

This line of questioning has its origins in \cite{JoPe98}, in which Jorgensen and Pedersen demonstrated Fourier bases for $L^2(\mu_{\lambda})$ when $\lambda = \frac{1}{2n}$ for each $n \in \mathbb{N}$.  They also showed that when $\lambda = \frac{1}{2n+1}$, there is no orthonormal basis (ONB) consisting of exponential functions, and in fact, every orthogonal collection of exponentials is finite when the denominator of $\lambda$ is odd.   The Fourier basis found in \cite[Cor.\! \!5.9]{JoPe98} is indexed by a discrete set $\Gamma \subset \mathbb{R}$ given by
\begin{equation}\label{Eqn:Gamma}
\Gamma = \Gamma\Bigl(\frac{1}{2n}\Bigr) = 
\Biggl\{ 
\sum_{i=0}^m a_i(2n)^i \: : \: a_i\in \Bigl\{0, \frac{n}{2}\Bigr\}, m \textrm{ finite }
\Biggr\}.
\end{equation}
For example, when $n=2$, i.e.  $\lambda = \frac{1}{4}$, the spectrum is \[ \Gamma\Bigl(\frac{1}{4}\Bigr) = 
\Biggl\{ 
\sum_{i=0}^m a_i 4^i \: : \: a_i\in \Bigl\{0,1 \Bigr\}
\Biggr\} = \{ 0, 1, 4, 5, 16, 17, 20, \ldots \}. \]
Note that the elements of $\Gamma$ are integers when $n$ is even and are all in $\frac{1}{2}\mathbb{Z}$ when $n$ is odd.  

\begin{remark}In the later sections, we will be considering other ONBs for the same Hilbert space $L^2(\mu)$.  To distinguish the original Jorgensen-Pedersen basis, we will call $\Gamma$ the \textit{canonical} ONB for $L^2(\mu)$ and label the others \textit{alternate} ONBs. \end{remark}

 We observe that $\Gamma$ has a self-similarity by scaling in the large.  For instance, $\Gamma$ is invariant under scaling by the value $2n$, the reciprocal of $\lambda$.  But there is even a stronger scaling invariance: \[ \Gamma = 2n\Gamma \sqcup \Big(2n\Gamma + \frac{n}{2}\Big).\]

Recent results have shown that if $\lambda > \frac12$, there is no Fourier basis \cite[Thm.\! \!3.4]{DHJ09}.  It has also been shown that if $\lambda = \frac{q}{2n}$ where $q$ is odd, the set of exponentials shown in \cite{JoPe98} to be an ONB for $\mu_{\frac{1}{2n}}$ are also orthogonal for $\mu_{\frac{q}{2n}}$ \cite[Thm.\! \!11.9]{JKS08}, \cite[Thm.\! \!1.2]{HuLa08}.  There are still a variety of open questions regarding the existence and classification of Fourier bases for Bernoulli convolution measures.  

We note that, since we wish to discuss ONBs for Bernoulli measures, we will be restricting here to the case where the scale factor $\lambda = \frac{1}{2n}$ for $n \in \mathbb{N}$.  To keep the notation simple, we will write $\mu$ for $\mu_{\frac{1}{2n}}$ and $X$ for $X_{\frac{1}{2n}}$.  We will use the notation $e_{\gamma}$ to represent the exponential function $e^{2\pi i \gamma (\cdot)}$.  If a collection of exponential functions  $\{e_{\gamma}\,:\, \gamma \in \Gamma\}$ indexed by the set $\Gamma$ forms an orthonormal basis for some $L^2(\mu_{\lambda})$, we call $\Gamma$ a \textit{spectrum} for $L^2(\mu_{\lambda})$.  Also, for ease of notation, we will write $E(\Gamma)$ for the set of exponentials $\{e_{\gamma}\,:\, \gamma \in \Gamma\}$.  

\subsection{The Fourier transform $\widehat{\mu}$ and its zero set}

Given the invariance equation \eqref{eqn:mu_invariance}, there is a standard convenient expression for the integral of an exponential function $e_t$.  We denote the resulting function in $t$ by $\widehat{\mu}$ since this produces a Fourier transform of $\mu$:

\begin{eqnarray*}  \widehat{\mu}(t) &=& \int_{X_{\lambda}} e^{2\pi i x t} \dd \mu(x) \\ &\underset{\eqref{eqn:mu_invariance}}{=}& \frac12 \int_{X_{\lambda}} \Big(e^{2\pi i \lambda (x+1) t} + e^{2\pi i \lambda (x-1) t} \Big) \dd \mu(x) \\ &=& \cos(2\pi \lambda t) \widehat{\mu}(\lambda t) \\ &=& \cos(2\pi \lambda t) \cos(2\pi \lambda^2 t) \widehat{\mu}(\lambda^2t)\\&=& \vdots\quad . \end{eqnarray*}
Continuing the iteration, we find an infinite product formula for $\widehat{\mu}$:
\begin{equation}\label{eqn:muhat}
\widehat{\mu}(t)=  \prod_{k=1}^{\infty} \cos(2\pi \lambda^kt).
\end{equation}
Given exponential functions $e_{\gamma}$ and $e_{\gamma'}$, we note that 
\begin{eqnarray}\label{eqn:innerproduct}\nonumber  \langle e_{\gamma'}, e_{\gamma} \rangle_{L^2(\mu)} &=& \int_{X} e^{2\pi i (\gamma - \gamma') x} \dd \mu(x) \\ &=& \widehat{\mu}(\gamma - \gamma'). \end{eqnarray}  We are considering orthogonal collections of exponential functions, so we consider the zeroes of the function $\widehat{\mu}$.  By Equation \eqref{eqn:muhat}, $\widehat{\mu}$ is zero if and only if one of the factors in the infinite product is zero.  (This is not difficult to show, but see \cite[Lem.\! 11.2]{JKS08} for explicit details.)   The cosine function is zero at the odd multiples of $\frac{\pi}{2}$, which yields the set of zeroes for $\widehat{\mu}$, denoted $\mathcal{Z}$:
\begin{equation}\label{Eqn:FourierZeroes}
\mathcal{Z}(\widehat{\mu}_{\frac{1}{2n}} )= \Biggl\{ \frac{(2n)^k(2m+1)}{4} \: \Big| \: m\in\Z, k \geq 1\Biggr\}.
\end{equation}

Given a discrete set $\Gamma$, then, the collection of exponential functions $E(\Gamma)$ is an orthogonal collection if $ \langle e_{\gamma'}, e_{\gamma} \rangle_{L^2{\mu})} = \delta_{\gamma, \gamma'}, $ which occurs if and only if \[ \gamma - \gamma' \in \mathcal{Z} \quad \textrm{ for all } \gamma, \gamma' \in \Gamma \textrm{ with } \gamma \neq \gamma'.\]  Note that the set from Equation \eqref{Eqn:Gamma} does indeed satisfy this condition.

We computed earlier that \[ \widehat{\mu}(t) = \cos(2\pi \lambda t) \widehat{\mu}(\lambda t).\]  Making use of our current situation for which $\lambda = \frac{1}{2n}$ and a change of variables, we also have the useful equation:
\begin{equation}\label{Eqn:MuHat2nt}
\widehat{\mu}(2nt) = \cos(2\pi t)\: \widehat{\mu}(t)
\end{equation}
for all $t\in\R$.  The utility of this formula comes when $t$ is an integer, as then we can actually cancel the factor $2n$, i.e. 
\[
\widehat{\mu}(2nt) = \widehat{\mu}(t).
\]

\section{Cuntz Operators }\label{sec:Cuntz}

In this section, we define operators on $L^2(\mu)$ which satisfy Cuntz relations.  The resulting Proposition \ref{Prop:decomp} will provide an orthogonal decomposition of the Hilbert space $L^2(\mu)$, which allows us to study Fourier bases for $L^2(\mu)$ in Section \ref{sec:U}.  In many cases, there may be more than one spectrum for $L^2(\mu)$, and if so, we can define quite naturally the unitary operator $U$ mapping one ONB to another.  Finally, we examine the matrix for $U$ with respect to the orthogonal decomposition given by the Cuntz operators.

\subsection{The Cuntz relations}

It is known that for every integer $N > 1$, there is a rich variety of non-equivalent representations of the Cuntz relations on $N$ generators \cite{Cun77}. The corresponding $C^*$-algebra on the $N$ generators is denoted $O_N$. Even the case $N = 2$ offers a rich family of such representations with applications to analysis of wavelets and, more generally, to fractals generated by two transformations; see Definition \ref{Defn:CuntzRel} below.   We note that the use of the Cuntz relations is not restricted to the harmonic analysis of \textit{affine} fractals; the relations can be applied much more generally.  Our purpose here is to isolate a particular representation of $O_2$ which lets us produce powerful algorithms for Fourier expansions in $L^2(\mu)$. 

 \begin{definition}\label{Defn:CuntzRel}  We say that operators $S_0, S_1$ on $L^2(\mu)$ satisfy \textit{Cuntz relations}
if 
\begin{enumerate}  
\item $S_0S_0^* + S_1S_1^* = I$,  
\item $S_i^*S_j = \delta_{i,j}I$ \;\; for $i,j = 0, 1$.
\end{enumerate}
When these relations hold, we have a representation of the Cuntz algebra $O_2$.
\end{definition}

We begin by defining scaling operators $S_0$ and $S_1$ on the Hilbert space $L^2(\mu)$ and showing that they satisfy Definition \ref{Defn:CuntzRel}; this representation yields an orthogonal decomposition of $L^2(\mu)$.  To motivate the connection between the harmonic analysis of affine fractals and operator representations within a $C^*$-algebra,  we look to the recurring notion of self-similarity.  The self-similarity or scaling invariance present in $\Gamma$ given by 
\[
\Gamma = 2n\Gamma \sqcup \Bigl(2n\Gamma +\frac{n}{2}\Bigr).
\]
has a representation via the Cuntz relations on bounded operators on $L^2(\mu)$.  The introduction of representations of the Cuntz relations allows us to transfer spectral duality questions to a framework in geometry of Hilbert space.  

Given the Bernoulli measure $\mu$ with scale factor $\frac{1}{2n}$, let $\Gamma$ be the canonical spectrum.      We define $S_0, S_1$ on the canonical ONB $E(\Gamma)$:  
\begin{equation}\label{Eqn:S0}
S_0e_{\gamma} = e_{2n\gamma}
\end{equation}
and
\begin{equation}\label{Eqn:S1}
S_1e_{\gamma} = e_{\frac{n}{2} + 2n \gamma}.
\end{equation}

Using Parseval's identity, we find that $S_0$ and $S_1$ are isometries since they each map $E(\Gamma)$ into a subset of itself.  The range of $S_0$ is the subspace spanned by $E(2n\Gamma)$, and the range of $S_1$ is the subspace spanned by $E(\frac{n}{2} + 2n\Gamma)$.  It is apparent that the range spaces are orthogonal and that $L^2(\mu)$ is the direct sum of these spaces.  We can extend this decomposition by observing that $S_0$ and $S_1$  satisfy the Cuntz relations from Definition \ref{Defn:CuntzRel}.

In order to verify that $S_0$ and $S_1$ verify the Cuntz relations, we first compute their adjoints.   Given $\xi, \gamma$ from $\Gamma$,
\begin{equation*}
\langle e_{\xi}, S_0^*e_{\gamma}\rangle = \langle S_0e_{\xi}, e_{\gamma}\rangle = \widehat{\mu}(2n\xi - \gamma). \end{equation*}
If $\gamma$ is not in $2n\Gamma$ then the exponentials are orthogonal.  If, on the other hand, $\gamma = 2n\gamma'$ for some $\gamma' \in \Gamma$, then by Equation (\ref{Eqn:MuHat2nt}),
\begin{equation*}
\widehat{\mu}(2n\xi - 2n\gamma') = \widehat{\mu}(\xi - \gamma') = \langle e_{\xi}, e_{\gamma'}\rangle.
\end{equation*}
Thus, we have \begin{equation}\label{Eqn:S_0^*}
S_0^*e_{\gamma} = \left\{ \begin{matrix} e_{\frac{\gamma}{2n}} & \textrm{when } \gamma \in 2n\Gamma \\
0 & \textrm{otherwise.} \end{matrix} \right.
\end{equation}

Similarly, 
\[ \langle e_{\xi}, S_1^*e_{\gamma} \rangle = \langle S_1e_{\xi}, e_{\gamma} \rangle = \widehat{\mu}\Bigl(2n\xi + \frac{n}{2} - \gamma\Bigr). \]
Here, we get zero when $\gamma \in 2n\Gamma$.  If $\gamma = 2n\gamma' + \frac{n}{2}$, 
\[ \widehat{\mu}\Bigl(2n\xi + \frac{n}{2} - 2n\gamma' - \frac{n}{2}\Bigr) = \widehat{\mu}(2n(\xi - \gamma')) = \widehat{\mu}(\xi - \gamma') = \langle e_{\xi}, e_{\gamma'} \rangle .\]
This gives
\begin{equation}\label{Eqn:S_1^*} 
S_1^*e_{\gamma} =  \left\{ \begin{matrix} e_{\frac{\gamma- \frac{n}{2}}{2n}} & \textrm{when } \gamma \in \frac{n}{2} + 2n\Gamma \\
0 & \textrm{otherwise.} \end{matrix} \right.
\end{equation}

\begin{proposition}\label{prop:cuntz}
The operators $S_0$ and $S_1$ defined in Equations (\ref{Eqn:S0}) and (\ref{Eqn:S1}) satisfy the Cuntz relations in Definition \ref{Defn:CuntzRel}.
\end{proposition}
\begin{proof} We start with the second part:  $S_i^*S_j = \delta_{i,j}I$ for $i,j = 0, 1$.  Then we show that $S_0S_0^* + S_1S_1^* = I$.

(2) \;  It follows directly from the computation of the adjoints $S_0^*$ and $S_1^*$ that applying either $S_0^*$ to an element in the range of $S_1$ or applying $S_1^*$ to an element in the range of $S_0$ yields $0$.  Also from the definitions, we see that the adjoints act as left inverses to the original operators, i.e. $S_0^*S_0 = I = S_1^*S_1$. 

(1) \;  Given the isometries $S_i$, $i=0,1$, $S_iS_i^*$ is a projection onto the range of $S_i$ since \[ (S_iS_i^*)(S_iS_i^*) = S_i(S_i^*S_i)S_i^* = S_iS_i^*\] by (2) above.  We have already shown that the ranges are orthogonal and their direct sum is all $L^2(\mu)$, hence $S_0S_0^* + S_1S_1^* = I$.
\end{proof}

\begin{remark}  Proposition \ref{prop:cuntz}  can be seen as a special case ($N=2$) of Theorem 2.3 in  \cite{DJ10}, which uses the Cuntz relations to formulate an algorithmic approach for constructing different Fourier ONBs for a fixed IFS measure $\mu$.     We note here, however, that the spectral scaling by $p$ which we introduce in Section \ref{sec:U} to build the operator-fractal $U$  is not part of a Cuntz algebra representation.    
\end{remark}

We will use $S_0$ and $S_1$ to describe an orthogonal decomposition of $L^2(\mu)$.   The first step will be to define the following decomposition of the canonical spectrum $\Gamma$.  Let 
\begin{eqnarray*}  \Gamma_0 &=& \Gamma \setminus 2n\Gamma = \frac{n}{2} + 2n\Gamma \\ \Gamma_1 &=& 2n\Gamma \setminus (2n)^2\Gamma = 2n\Bigl(\frac{n}{2} + 2n\Gamma\Bigr) \\ \vdots & & \vdots \\ \Gamma_k &=& (2n)^k\Gamma \setminus (2n)^{k+1}\Gamma = (2n)^k\Bigl(\frac{n}{2} + 2n\Gamma\Bigr).
\end{eqnarray*} 
It is readily verified that the sets $\Gamma_k$ are disjoint and form a partition of $\Gamma \setminus \{0\}$.  

The elements of $\Gamma$ are in a natural one-to-one correspondence with infinite sequences of bits from $\{0, \frac{n}{2}\}$ having finitely many nonzero elements:
\begin{equation}\label{Eqn:SeqCor}
 \gamma = \sum_{i=0}^m a_i(2n)^i \leftrightarrow (a_0, a_1, a_2, \ldots, a_m, 0, 0, 0, \ldots),
\end{equation}
where $a_i\in \{0, \frac{n}{2}\}$.  We can describe $\Gamma_k$  as the set of all elements from $\Gamma$ which correspond to a sequence having exactly $k$ leading zeroes (note the index starts at $0$): \[ (0, 0, \ldots, 0, \frac{n}{2}, a_{k+1}, \ldots, a_m, 0, 0, 0, \ldots ). \]

The operators $S_0$ and $S_1$ can also be associated to operators on sequences via their behavior on the ONB elements $e_{\gamma}$.  The operator $S_0$ takes $e_{\gamma}$ to $e_{2n\gamma}$, which corresponds to a right shift on the sequence associated to  $\gamma$.   Since $S_1$ takes $e_{\gamma}$ to $e_{\frac{n}{2} + 2n\gamma}$, we associate the action of $S_1$ to the operator that shifts the sequence for $\gamma$ right and places $\frac{n}{2}$ in the vacant position.  
\begin{eqnarray*} \gamma &\leftrightarrow& (a_0, a_1, a_2, \ldots, a_m, 0, 0, 0, \ldots) \\ 2n\gamma & \leftrightarrow & (0, a_0, a_1, \ldots, a_m, 0, 0, \ldots) \\ \frac{n}{2} + 2n\gamma & \leftrightarrow & (\frac{n}{2}, a_0, a_1, \ldots, a_m, 0, 0, \ldots) \end{eqnarray*}

We combine the associations above to prove the following.

\begin{lemma}\label{lem:Wk}  The closed span of $E(\Gamma_k)$ is the subspace $S_0^kS_1(L^2(\mu))$, for each $k=0, 1, 2, \ldots$. 
\end{lemma}

\begin{proof} We show that  $e_{\xi} \in E(\Gamma_k)$ if and only if  $e_{\xi} = S_0^kS_1e_{\gamma}$ for some $\gamma \in \Gamma$.  

From the definition,  $\gamma' \in \Gamma_k$ if and only if $\gamma' = (2n)^k(\frac{n}{2} + 2n\gamma)$ for some $\gamma \in \Gamma$.  Then, for this $\gamma$,  \[S_0^kS_1e_{\gamma} = S_0^ke_{\frac{n}{2} + 2n\gamma} = e_{(2n)^k(\frac{n}{2} + 2n\gamma)} = e_{\gamma'}.\]  Thus, \[ E(\Gamma_k) = S_0^kS_1(E(\Gamma)), \quad k=0, 1, 2, \ldots .\]  Since $E(\Gamma)$ is an orthonormal basis for $L^2(\mu)$, the closed span of $E(\Gamma_k)$ is the subspace $S_0^kS_1(L^2(\mu))$. \end{proof}

For $k = 0, 1, 2, \ldots$, let $W_k$ denote the following set: 
\[ W_k:= S_0^kS_1(L^2(\mu)).\]
Note that $W_{k+1} = S_0 W_k$.  By the previous lemma, each $W_k$ is spanned by $E(\Gamma_k)$.   Since the sets $\Gamma_k$ form a partition of $\Gamma \setminus \{0\}$, the $W_k$ subspaces are mutually orthogonal, and their direct sum yields $L^2(\mu) \ominus sp\{e_0\}$.  We record the orthogonal decomposition of $L^2(\mu)$ in the following proposition.

\begin{proposition}\label{Prop:decomp}
Given the spaces \begin{equation}\label{eqn:Wk} W_k = S_0^kS_1(L^2(\mu)) \textrm{  for  } k=0, 1, 2, \ldots, \end{equation} we have \begin{equation} L^2(\mu) = sp\{e_0\} \oplus \bigoplus_{k=0}^{\infty} W_k. \end{equation}
\end{proposition}

\section{The operator-fractal $U$}\label{sec:U}

There are many cases (for example, \cite{LaWa02, DJ09, DHS09, JKS11a}) in which a fractal $L^2(\mu)$ Hilbert space has more Fourier bases than just the canonical one.  We will refer to these as \textit{alternate ONBs} for the space.  One technique in the papers above for finding alternate ONBs is to scale the spectrum $\Gamma$ of the canonical ONB by an odd integer $p$.  The new set $p\Gamma$ always yields an orthogonal collection of exponentials and often yields an alternate ONB.  Examples include $p=5$ for the $\lambda = \frac14$ measure \cite[Prop.\! 5.1]{DJ09} and $p=3$ for the $\lambda = \frac18$ measure \cite[Thm.\! 3.5]{JKS11a}.

In this section, we will use pairs of canonical and alternate ONBs to define the unitary operator $U$.  We will demonstrate that $U$ has inherent fractal-like properties, and we will find additional structure in the $\lambda = \frac14$, $p=5$ example.  

Consider an arbitrary pair for $\lambda = \frac{1}{2n}$ where the canonical set $\Gamma = \Gamma(\frac{1}{2n})$ and the scaled set $p\Gamma$ are both spectra.  Define the operator $U$ on  $\Gamma$ by \begin{equation}\label{Eqn:U}  Ue_{\gamma} = e_{p\gamma}.\end{equation}

Since $U$ maps an ONB to another ONB, it is a unitary operator on $L^2(\mu)$.  Note that $U^2(e_{\gamma})\neq e_{p^2\gamma}$ in general.  Computation of powers of $U$ requires an infinite expansion with respect to the $E(\Gamma)$ ONB, i.e.   $U^2 e_{\gamma} = \sum_{\xi \in \Gamma} \langle e_{\xi}, Ue_{\gamma} \rangle Ue_{\xi}$.  Toward this end, it becomes useful to study the matrix representation for $U$ with respect to $E(\Gamma)$.  We find a surprising self-similarity in the matrix for $U$.

Given $\gamma, \xi \in \Gamma$, the matrix entries for the operator $U$ are
 \[
 U(\xi, \gamma) = \langle e_{\xi}, Ue_{\gamma} \rangle = \widehat{\mu}(p\gamma - \xi).
 \]  If we reorder the index set $\Gamma$ by the subsets $\Gamma_k$, we find that the matrix for $U$ has a block structure.  The block structure arises because the subspaces $W_k$ and $U(W_{k'})$ are orthogonal for $k \neq k'$, but we will in fact prove the stronger result that each space $W_k$ is invariant under $U$.  

\subsection{Invariance properties of $U$}
In this section, we examine the interplay of $U$ and the Cuntz subspaces.  We find that the subspaces are invariant under $U$.
From now on, we use $\mathcal{L}$ to denote 
\[\mathcal{L} = L^2(\mu) \ominus sp\{e_0\} = \bigoplus_{k=0}^{\infty} W_k.\]  
\begin{lemma}\label{lem:Ukinvariant}For each $k \geq 0$, $S_0^k(\mathcal{L})$ is invariant under $U$.
\end{lemma}
\begin{proof}$ $\\
\textbf{Case 1:  }$k = 0$.  We show that $Ue_{\gamma}$ is in the orthogonal complement of $sp\{e_0\}$ if and only if $\gamma \neq 0$.  This guarantees that $Ue_{\gamma} \in \mathcal{L}$.  

If $\gamma = 0$ then $Ue_{\gamma} = e_0$.  Conversely, if $\gamma \neq 0$, then $Ue_{\gamma} = e_{p\gamma}$ where $p\gamma$ is a nonzero element of $p\Gamma$.  Then, since $E(p\Gamma)$ is an ONB which contains $e_0$, we have $\langle e_{p\gamma}, e_0 \rangle = 0$. 

\noindent\textbf{Case 2:  }$k \geq 1$. Let $\gamma \in \Gamma \setminus \{0\}$.   By the definition of $S_0$ in \eqref{Eqn:S0} and by Equation \eqref{eqn:innerproduct}, 
 \begin{eqnarray*} US_0^ke_{\gamma} &=& e_{p (2n)^k \gamma}\\ &=& \sum_{\xi' \in \Gamma \setminus \{0\}} \widehat{\mu}(p (2n)^k \gamma - \xi') e_{\xi'}. \end{eqnarray*}
In the expansion above, we leave out $\xi'=0$ since $US_0^k\mathcal{L} \subseteq US_0\mathcal{L} \subseteq \mathcal{L}$.

Consider the terms in the sum for which $\xi' \notin (2n)^k\Gamma$.  Recalling that the sets $\Gamma_i, i \in \mathbb{N}_0$ partition $\Gamma \setminus \{0\}$, we find that $\xi' \notin (2n)^k\Gamma$ is equivalent to $\xi'$ being an element from $\Gamma_{i_0}$ for some $i_0 < k$.  Then $\xi' = (2n)^{i_0}(\frac{n}{2} + 2n \xi)$ for some $\xi \in \Gamma$, and 
\begin{eqnarray*} p (2n)^k \gamma - \xi' 
&=& 
\frac{1}{4} \Bigr(p (2n)^{k} 4\gamma - 4 (2n)^{i_0} \Big(\frac{n}{2} + 2n\xi \Big) \Bigr)\\
&=&
\frac{1}{4}(2n)^{i_0} \Bigr(p (2n)^{k-i_0} 4\gamma - 2n - 8n\xi \Bigr)\\ 
&=& 
\frac{1}{4} (2n)^{i_0+1}\Bigr( p (2n)^{k-i_0-1} 4\gamma - 1- 4 \xi \Bigr) \\
&\in& 
 \mathcal{Z}, \textrm{ the $\widehat{\mu}$ zero set from \eqref{Eqn:FourierZeroes}.} \end{eqnarray*}
Recall that, whether $n$ is even or odd, $4\gamma$ and $4\xi$ are even integers.  Therefore, the quantity in parentheses in the last expression of the equation above is an odd integer.  

This proves that $\widehat{\mu}(p(2n)^k\gamma - \xi') = 0$ when $\xi' \notin (2n)^k\Gamma$, hence $US_0^ke_{\gamma'}$ is in the closed span of $E((2n)^k(\Gamma\setminus\{0\})$, which is exactly the space $S_0^k(\mathcal{L})$.  Therefore $S_0^k(\mathcal{L})$ is invariant under $U$.
\end{proof}

Next, we show that the $W_k$ subspaces are invariant under $U$.  

 \begin{proposition}\label{prop:Winvariant} Each space $W_k = S_0^kS_1(L^2(\mu))$ is invariant under $U$.
 \end{proposition} 
 \begin{proof}
 Let $e_{\gamma'} \in W_k$.  This implies that $e_{\gamma'} \in S_0^k(\mathcal{L}) \subseteq \mathcal{L}$.  We showed in Lemma \ref{lem:Wk} that $W_k = \overline{sp}\{E(\Gamma_k)\}$, so we have $\gamma' \in \Gamma_k$ is of the form $\gamma' = (2n)^k( \frac{n}{2} +2n\gamma)$ for some $\gamma \in \Gamma$.  Then,
 \[ Ue_{\gamma'} = e_{p\gamma'} = \sum_{\xi' \in \Gamma \setminus \{0\}} \widehat{\mu}(p\gamma'-\xi') e_{\xi'}, \] where, by Lemma \ref{lem:Ukinvariant}  the coefficients in this sum are possibly nonzero only for nonzero $\xi' \in (2n)^k\Gamma$.  We will show that the coefficients are also zero for $\xi' \in (2n)^{k+1}\Gamma \subset (2n)^k\Gamma$, hence the only possibly nonzero coefficients are from $\xi' \in \Gamma_k$.   
 
 Let $\xi' \in (2n)^{k+1}\Gamma$, i.e. $\xi' = (2n)^{k+1}\xi$ for some $\xi \in \Gamma$.  Then, $e_{\xi'} \in S_0^{k+1}(\mathcal{L})$ and
 
 \[ \widehat{\mu}(p\gamma'-\xi') = \widehat{\mu}\Big(p (2n)^k\Big(\frac{n}{2}+ 2n\gamma\Big) - (2n)^{k+1} \xi \Big).\]
We then compute that 
\[ 
p (2n)^k\Bigl(\frac{n}{2}+ 2n\gamma\Bigr) - (2n)^{k+1} \xi = \frac{1}{4} (2n)^{k+1}(p + 4p\gamma - 4\xi) \in \mathcal{Z}.
\]
Since $p\gamma' - \xi'$ is in the zero set, $Ue_{\gamma'}$ is orthogonal to each $e_{\xi'}$ where  $\xi' \in (2n)^{k+1}\Gamma$.  Thus, we have $Ue_{\gamma'} \in W_k$, so $W_k$ is invariant under $U$.
  \end{proof}

\subsection{Commutation relations for $U$ and $S_0$}\label{Subsec:Commute}

We can actually say something stronger than invariance about the interaction between $U$ and $S_0$.  When $n$ is even, the operators commute.  When $n$ is odd, they don't quite commute but do have a straightforward relationship.  We note that $U$ and $S_1$ do not commute.

\begin{theorem}\label{Prop:EvenCommute}
If $n$ is even, then the operators $S_0$ and $U$ commute on $L^2(\mu)$.
\end{theorem}

\begin{proof}  Let $\gamma\in\Gamma$.  Note first that $2n\gamma\in \Gamma$, while $p\gamma$ is not necessarily an element of $\Gamma$.   We must therefore expand both expressions $US_0 e_\gamma$ and $S_0Ue_{\gamma}$ in terms of the ONB $E(\Gamma)$.  First,
\begin{equation}\label{Eqn:S_0U}
\begin{split}
S_0 Ue_{\gamma} & = S_0 e_{p\gamma}
= S_0\Bigl( \sum_{\eta\in\Gamma} \langle e_{p \gamma}, e_{\eta}\rangle e_{\eta}\Bigr)
 =  \sum_{\eta\in\Gamma} \langle e_{p \gamma}, e_{\eta}\rangle e_{2n \eta}\\
& =  \sum_{\eta\in\Gamma} \widehat{\mu}(p \gamma - \eta) e_{2n \eta}.
\end{split}
\end{equation}
On the other hand,
\begin{equation}\label{Eqn:US_0_1}
US_0e_{\gamma}  = U e_{2n\gamma} = e_{p 2n \gamma}.
\end{equation}
Since $\Gamma = 2n\Gamma \sqcup \Bigl(\frac{n}{2}+2n\Gamma\Bigr) = 2n\Gamma \sqcup  \Gamma_0$, we can split the sum over $\Gamma$ into two terms:
\begin{equation}
\begin{split}
US_0e_{\gamma} 
& = \sum_{\xi\in\Gamma} \widehat{\mu}(p 2n \gamma - \xi)e_{\xi}\\
& =  \sum_{\xi\in 2n\Gamma} \widehat{\mu}(p  2n \gamma - \xi)e_{\xi} +  \sum_{\xi\in\Gamma_0} \widehat{\mu}(p  2n \gamma - \xi)e_{\xi}.
\end{split}
\end{equation}
However, by Lemma \ref{lem:Ukinvariant}, $\widehat{\mu}(p 2n \gamma - \xi) = 0$ when $\xi\in\Gamma_0$, hence \begin{equation}\label{Eqn:US_0_2}
\begin{split}
US_0e_{\gamma} 
& =  \sum_{\xi\in 2n\Gamma} \widehat{\mu}(p 2n \gamma - \xi)\:e_{\xi}\\
& =  \sum_{\xi'\in \Gamma} \widehat{\mu}(p 2n \gamma - 2n \xi')\:e_{2n \xi'}.\\
\end{split}
\end{equation}
Since $n$ is even, $\gamma\in \Z$ and $\xi\in \Z$.  By applying Equation (\ref{Eqn:MuHat2nt}), we have
\begin{equation}
US_0e_{\gamma} = \sum_{\xi'\in\Gamma} \widehat{\mu}(p \gamma - \xi')\:e_{2n \xi'},
\end{equation}
which agrees exactly with the expression for $S_0Ue_{\gamma}$ in Equation (\ref{Eqn:S_0U}).

Note that $U$ and $S_0$ commute when restricted to $\mathcal{L}$ as well.  
\end{proof}

In the previous theorem, we depended on the fact that $n$ was even to apply Equation \eqref{Eqn:MuHat2nt} in the last step in Equation (\ref{Eqn:US_0_2}).  In fact, if $n$ is not even, $U$ and $S_0$ do not quite commute, as the next theorem shows.  First, we point out some straightforward observations.
\begin{itemize}
\item \; When $n$ is odd and $\gamma\in \Gamma_0 = \frac{n}{2} + 2n\Gamma$, we see that $\gamma$ is \textit{not} an integer, and is in fact an element of $\frac12(2\N-1)$.  Therefore, $2n\gamma$ is an odd integer for $\gamma \in \Gamma_0$ and an even integer for $\gamma \in \Gamma \setminus \Gamma_0 = 2n\Gamma$.
\item\;   Let $n$ and  $p$ be odd and let $\gamma\in \Gamma_0$ and $\xi\in\Gamma$.  Then, since $2n\gamma \in \N$, 
\[ p\gamma - 2n \xi \in \frac{1}{2}(2\Z +1).\]   
\item \;  Given $n$ and $p$ odd, \[ \cos(2 \pi (p \gamma - \xi)) = \begin{cases} 1 & \gamma, \xi \in 2n\Gamma \\ 1 & \gamma, \xi \in \Gamma_0 \\ -1 & \gamma \in 2n\Gamma,\, \xi \in \Gamma_0 \\ -1 & \gamma \in \Gamma_0,\,  \xi \in 2n\Gamma. \end{cases} \]
\end{itemize}

\begin{theorem}\label{Prop:OddCommute}  Suppose $n$ is odd and $\gamma\in\Gamma$.  Write
\[ S_0 Ue_{\gamma} = s + w,\]
where $s\in S_0^2(L^2(\mu))$ and $w\in W_1 = S_0S_1(L^2(\mu))$.  Then
\begin{equation}
US_0e_{\gamma} = 
\begin{cases} s - w &  \textrm{ if } \gamma\in 2n\Gamma\\
-s + w & \textrm { if } \gamma\in \frac{n}{2} + 2n\Gamma = \Gamma_0.
\end{cases}
\end{equation}
\end{theorem}
\begin{proof} By Equation (\ref{Eqn:S_0U}) in Theorem \ref{Prop:EvenCommute}, we have
\begin{equation}
S_0 Ue_{\gamma} 
=  \sum_{\eta\in\Gamma} \widehat{\mu}(p \gamma - \eta) e_{2n \eta}
\end{equation}
Since  $2n\Gamma =  (2n)^2\Gamma \sqcup \Gamma_1$, we split up the sum into these two pieces, which we call respectively $s$ and $w$:
\begin{equation}\label{Eqn:S_0Uodd}
\begin{split}
S_0 Ue_{\gamma} 
& = \sum_{\eta\in\Gamma} \widehat{\mu}\Bigl(p\cdot\gamma - 2n\cdot \eta\Bigr) e_{(2n)^2\eta}\\
& \phantom{{===}} + \sum_{\eta\in\Gamma} \widehat{\mu}\Bigl(p\cdot \gamma - \frac{n}{2} - 2n\cdot \eta\Bigr) e_{2n(\frac{n}{2} + 2n\eta)}\\
& = s \oplus w \in S_0^2(L^2(\mu)) \oplus W_1.
\end{split}
\end{equation}

We split the sum for $US_0e_{\gamma}$ in the same fashion, yielding $\widetilde{s}$ and $\widetilde{w}$:
\begin{equation}\label{Eqn:US_0_3}
\begin{split}
US_0e_{\gamma} 
& = \sum_{\xi\in 2n\Gamma}  \widehat{\mu}(p  2n  \gamma - \xi)\:e_{\xi}\\
& = \sum_{\xi\in\Gamma}\widehat{\mu}(p  2n  \gamma - (2n)^2\xi)\:e_{(2n)^2\xi}\\
& \phantom{{===}}  + \sum_{\xi\in\Gamma}\widehat{\mu}\Bigl(p  2n  \gamma - (2n)\Bigl(\frac{n}{2} + 2n \xi\Bigr)\Bigr)\:e_{2n(\frac{n}{2} + 2n\cdot \xi)}\\ &= \widetilde{s} \oplus \widetilde{w} \in S_0^2(L^2(\mu)) \oplus W_1\\
\end{split}
\end{equation}

\noindent \textbf{Case 1:  }$\gamma\in 2n\Gamma$.  Let $\gamma = 2n\gamma' \in 2n\Gamma$.  Then, making use of Equation \eqref{Eqn:MuHat2nt}, we see that the coefficients for $\widetilde{s}$ in Equation \eqref{Eqn:US_0_3} are
\begin{equation*}
\begin{split}
&\widehat{\mu}(p (2n)^2 \gamma' - (2n)^2 \xi) \\
&=\cos(2\pi (p 2n\gamma' - 2n\xi)) \:\widehat{\mu}(p 2n\gamma' - 2n\xi) \\ 
&=
\widehat{\mu}(p  2n\gamma' - 2n\xi),
\end{split}
\end{equation*}
since $p 2n\gamma' - 2n\xi \in \Z$.  These are the same as the coefficients for $s$ from Equation \eqref{Eqn:S_0Uodd}, so $\widetilde{s} = s$.

The coefficients for $\widetilde{w}$ in Equation \eqref{Eqn:US_0_3} are
\begin{equation*}
\begin{split}
& \widehat{\mu}\Bigl(p (2n)^2 \gamma' - 2n\Bigl(\frac{n}{2} + 2n\xi\Bigr)\Bigr) \\
& =
\cos\Bigl(2\pi\Bigl(p 2n\gamma' - \frac{n}{2} - 2n\xi\Bigr)\Bigr)\:\widehat{\mu}\Bigl(p 2n\gamma' - \frac{n}{2} - 2n\xi\Bigr) \\ 
& = 
(-1) \widehat{\mu}\Bigl(p  2n\gamma' - \frac{n}{2} - 2n\xi\Bigr). 
\end{split}
\end{equation*}
The cosine factors are all $-1$ since $p 2n\gamma' - \frac{n}{2} - 2n\xi$ is an element of $\frac12(2\Z+1)$.  Therefore, the coefficients for $\widetilde{w}$ are exactly the opposite sign from those for $w$, i.e. $\widetilde{w} = -w$.

\noindent \textbf{Case 2:  } $\gamma\in \Gamma_0 = \frac{n}{2} + 2n\Gamma$

Let $\gamma = \frac{n}{2} + 2n\gamma'$.  Again using Equation \eqref{Eqn:MuHat2nt}, we find that the coefficients for $\widetilde{s}$ are
\begin{equation*}
\begin{split} 
& \widehat{\mu}\Bigl(p2n\Bigl(\frac{n}{2} + 2n\gamma'\Bigr) - (2n)^2\xi\Bigr) \\
&=
\cos\Bigl(2\pi\Bigl( \frac{n}{2} + 2n\gamma' - 2n\xi\Bigr)\Bigr) \:\widehat{\mu}\Bigl( \frac{n}{2} + 2n\gamma' - 2n\xi\Bigr) \\ 
&= 
(-1)\widehat{\mu}\Bigl( \frac{n}{2} + 2n\gamma' - 2n\xi\Bigr)
\end{split}
\end{equation*}
The coefficients for $\widetilde{s}$ are the opposite sign to those for $s$, so $\widetilde{s} = -s$.  

The coefficients for $\widetilde{w}$ are 
\begin{equation*}
\begin{split} 
&
\widehat{\mu}\Bigl(p 2n\Bigl(\frac{n}{2} + 2n\gamma'\Bigr) - 2n\Bigl(\frac{n}{2} + 2n\xi\Bigr)\Bigr) \\
&=
\widehat{\mu}\Bigl(2n\Bigl(\frac{p-1}{2}n + p2n\gamma' - 2n\xi\Bigr)\Bigr) \\ 
&=
\widehat{\mu} \Bigl(\frac{p-1}{2}n + p2n\gamma' - 2n\xi\Bigr) \\ 
&= 
\widehat{\mu}\Bigl(p \gamma - \frac{n}{2} - 2n\xi\Bigr)
\end{split}
\end{equation*}
since $\frac{p-1}{2}$ is an integer.   These are the same coefficients as those for $w$, so $\widetilde{w} = w$.
\end{proof}

\begin{corollary}\label{cor:commutant}
When $n$ is odd, the commutator $US_0-S_0U$ is given by \[(US_0 - S_0U)  e_{\gamma} = \begin{cases} -2w & \gamma \in 2n\Gamma \\ -2s & \gamma \in \Gamma_0 \end{cases}, \]  where $s$ and $w$ are as defined in Proposition \ref{Prop:OddCommute}. \end{corollary}


The isometry operators $S_0, S_1$ have range projections $S_0S_0^*$ and $S_1S_1^*$ respectively.  The preceding results lead us to find that $U$ commutes with the range projections.

\begin{corollary}
The unitary operator $U$ commutes with $S_0S_0^*$ and $S_1S_1^*$.
\end{corollary}

\begin{proof}
This follows immediately because the spaces $W_0$ and $S_0(L^2(\mu))$ are invariant under $U$.
\end{proof}

\subsection{The matrix for $U$}\label{Subsec:Matrix}
We showed in Proposition \ref{Prop:decomp} that the $W_k$ subspaces are orthogonal, and since Proposition \ref{prop:Winvariant} showed that each $W_k$ is preserved by the operator $U$, we find that the infinite matrix for $U$, when $\Gamma$ is ordered by the $\Gamma_k$ sets, has a block matrix structure.  We note that, with the exception of the upper left corner, all of the blocks are infinite.

\begin{equation}\label{Eqn:MatrixU_02}
U \sim \begin{tabular}{c||c|c|c|c|c|c}
                         &\:\: $0$ \:\:& $\Gamma_0$  & $\Gamma_1$ & $\Gamma_2$ & $\Gamma_3$ & $\cdots$\\
\hline
\hline
$0$ & $1$ & $0$ & $0$ & $0$ & $0$& $\cdots$\\ 
\hline $\Gamma_0$  &$0$ & $U|_{W_0}$ &    $0$          &       $0$    &     $0$    &    $\cdots$ \\
\hline 
$\Gamma_1$  & $0$ &  $0$           &$U|_{W_1}$ &       $0$    &     $0$    &    $\cdots$ \\
\hline
$\Gamma_2$  & $0$ &  $0$           &   $0$            &        $U|_{W_2}$   &     $0$    &    $\cdots$ \\
\hline
$\Gamma_3$  & $0$ &  $0$           &   $0$            &        $0$                &    $U|_{W_3}$  &    $\cdots$ \\
\hline
                        & \vdots &  $\vdots$   &    $\vdots$   &      $\vdots$           &   $\vdots$         &  $\ddots$\\
\end{tabular}
\end{equation}

Next, we show that, in fact, the infinite blocks above are all the same.  For each $k \in \mathbb{N}_0$, consider the matrix of the operator $U|_{W_k}$ with respect to the indexing set $\Gamma_k$.  Let $\xi',\gamma' \in \Gamma_k$, so they have the form
\[\xi' =(2n)^k\Bigl( \frac{n}{2} + 2n\xi\Bigr)\]
and
\[\gamma' = (2n)^k\Bigl(\frac{n}{2} + 2n\gamma\Bigr),\] 
where $\xi, \gamma\in\Gamma$.  In the next lemma, we show that each $(\xi',\gamma')$ entry of $U|_{W_k}$ is equal to the corresponding $(\xi',\gamma')$ entry for $U|_{W_0}$.  In other words, the infinite blocks in the matrix for $U$ are all the same.

\begin{lemma}For all $k \in \mathbb{N}_0$, the infinite matrix for $U|_{W_k}$ has the same entries as the matrix for $U|_{W_{0}}$.
\end{lemma}
\begin{proof}
  Define $\xi', \gamma'$ as above.  The $(\xi', \gamma')$ entry of the  matrix for $U|_{W_k}$ is
\begin{equation*}
U|_{W_k}(\xi', \gamma')
= \Bigl\langle e_{(2n)^k(n/2 +2n\xi)}, Ue_{(2n)^k(n/2 +2n\gamma)}\Bigr\rangle_{L^2(\mu)}.
\end{equation*}
We now compute $U|_{W_k}(\xi', \gamma')$:
\begin{equation}\label{Eqn:UMatrixInnerProduct}
\begin{split}
& \Bigl\langle e_{(2n)^k(n/2 +2n\xi)}, Ue_{(2n)^k(n/2 +2n\gamma)} \Bigr\rangle_{L^2(\mu_)}\\
& =
\widehat{\mu}\Biggl( p\cdot (2n)^k\Bigl(\frac{n}{2} +2n\gamma\Bigr)-(2n)^k\Bigl(\frac{n}{2} +2n\xi\Bigr)  \Biggr) \\
& =
\widehat{\mu}\Biggl( (2n)^k \Biggl( \frac{(p-1)n}{2} + 2np\gamma - 2n\xi \Biggr) \Biggr).
\end{split}
\end{equation}
By the invariance property of $\mu$ in \eqref{Eqn:MuHat2nt}, 
\[
\widehat{\mu}(2nt) = \cos(2\pi t)\widehat{\mu}(t).
\]
If $t$ is an integer, then $\widehat{\mu}(2nt) = \widehat{\mu}(t)$.  Since $p-1$ is even, the number $t$ given by
\[
t = \frac{(p-1)n}{2} + 2np\gamma - 2n\xi
\] in Equation (\ref{Eqn:UMatrixInnerProduct}) is a whole number,  and thus we can disregard the $k$ powers of $2n$.  Therefore, we reach our desired conclusion:  
\begin{eqnarray*}
U|_{W_k}(\xi',\gamma')
& =& \widehat{\mu}\Biggl( (2n)^k \Biggl( \frac{(p-1)n}{2} + 2np\gamma - 2n\xi \Biggr) \Biggr)\\
& =& \widehat{\mu}\Biggl(\frac{(p-1)n}{2} + 2np\gamma - 2n\xi \Biggr) \\
& =& U|_{W_0}(\xi', \gamma'),
\end{eqnarray*}
since, in $W_0$, $\xi' = \frac{n}{2} + 2n\xi$ and $\gamma' = \frac{n}{2} + 2n\gamma$.
\end{proof}
Since the spaces $W_k$ are mutually orthogonal and are each invariant under the action of $U$, the infinite matrix representation for $U$ has the following block form.  As before, the index sets have been written across the top row and left column.
 
 \begin{theorem}\label{thm:block}  The matrix representation of the operator $U$, for the given ordering of the index set $\Gamma$, has a block diagonal structure of the form 
\begin{equation}\label{Eqn:MatrixU_03}
U \simeq  \begin{tabular}{c||c|c|c|c|c|c}
                         & \:\:$0$ \:\:& $\Gamma_0$  & $\Gamma_1$ & $\Gamma_2$ & $\Gamma_3$ & $\cdots$\\
\hline
\hline
$0$ & $1$ & $0$ & $0$ & $0$ & $0$& $\cdots$\\ 
\hline $\Gamma_0$  &$0$ & $U|_{W_0}$ &    $0$          &       $0$    &     $0$    &    $\cdots$ \\
\hline 
$\Gamma_1$  & $0$ &  $0$           &$U|_{W_0}$ &       $0$    &     $0$    &    $\cdots$ \\
\hline
$\Gamma_2$  & $0$ &  $0$           &   $0$            &        $U|_{W_0}$   &     $0$    &    $\cdots$ \\
\hline
$\Gamma_3$  & $0$ &  $0$           &   $0$            &        $0$                &    $U|_{W_0}$  &    $\cdots$ \\
\hline
                        & \vdots &  $\vdots$   &    $\vdots$   &      $\vdots$           &   $\vdots$         &  $\ddots$ \\ 
\end{tabular} \quad . \end{equation}
In other words, \[ U \simeq P_{e_0} \oplus \bigoplus_{k=0}^{\infty} U|_{W_0},\] where $P_{e_0}$ is the orthogonal projection onto $sp\{e_0\}$.
\end{theorem}

\begin{remark}  Making use again of Equation \eqref{Eqn:MuHat2nt}, and depending on the values of $n$ and $p$, we can write the entries of the matrix blocks as  
\begin{eqnarray*} \widehat{\mu}\Biggl(\frac{(p-1)n}{2} + 2np\gamma - 2n\xi \Biggr) &=& \widehat{\mu}\Biggl( (2n)\Bigl( \frac{p-1}{4} + p\gamma - \xi \Bigr) \Biggr)\\ &=& \pm \widehat{\mu}\Biggl( \frac{p-1}{4} + p\gamma - \xi \Biggr).\end{eqnarray*}
\end{remark}

\subsection{The individual blocks in the $p=5, n=2$ case}\label{Subsec:SpecialMatrix}

In this section, we consider the matrix of $U$ in the special case $p=5$ and $n=2$.  Here, $U(e_{\gamma}) = e_{5\gamma}$, and we zoom in on the subspace $W_0 = S_1(L^2(\mu_{\frac14}))$.  We will study the structure of one of the blocks in the matrix form for $U$ in Equation \eqref{Eqn:MatrixU_03}.  We find that, with a specified ordering of the index set of a block defined by the Cuntz operators $S_0$ and $S_1$, there is again a self-similar block structure in the matrix of $U|_{W_0}$.  

When $n=2$, the zero set for $\widehat{\mu}$ is 
\[
 \mathcal{Z}(\mu_{\frac{1}{4}})
 = 
 \Bigl\{ 4^k \cdot (2\ell + 1)\,:\, \ell \in \mathbb{Z}, k \in \mathbb{N}_0 \Bigr\}.
 \]
 The canonical ONB is 
 \[ 
 \Gamma = 
\Biggl\{ 
\sum_{i=0}^m a_i 4^i \: : \: a_i\in \{0, 1\}
\Biggr\} = \{0, 1, 4, 5, 16, 17, 20, 21, \ldots\}.
\]
By \eqref{Eqn:UMatrixInnerProduct}, the $(\xi',\gamma')$ entries in $U|_{W_k}$ are given by 
\begin{equation}\label{eqn:plusone}
U|_{W_k}(\xi', \gamma')
=
\widehat{\mu}_{\frac14}(1+5\gamma-\xi), 
\end{equation} 
where $\xi' = 4^k(1+4\xi)$ and $\gamma' = 4^k(1+4\gamma)$.   Each such block is indexed by the associated $\Gamma_k$, and by Theorem \ref{thm:block}, the entries do not depend on $k$.   It is sufficient, then, to examine the first block, indexed by $\Gamma_0$.  

We recall that the entries for the matrix of $U$ are $U(\xi, \gamma) = \widehat{\mu}(5\gamma - \xi)$.  This yields the expression of the restriction $U|_{W_0}$ as unitarily equivalent to the product of $U$ and a multiplication operator, which brings us to our surprising operator-fractal property of $U$.

\begin{theorem}\label{Thm:Multiplication}  Let $M_{e_1}$ be the multiplication operator by the function $e_1$ on $L^2(\mu)$, and let $P_{e_0}$ be the orthogonal projection onto $sp\{e_0\}$.  Then $U$ has the form
\begin{equation}\label{eqn:multiplicaton} 
U \simeq  P_{e_0} \oplus \bigoplus_{k=0}^{\infty} M_{e_1}U.
\end{equation} 
\end{theorem}\begin{proof}  First, observe that \[ U|_{W_0} = S_1^*US_1,\] since the range of the isometry $S_1$ is exactly the subspace $W_0$.  We see that
\begin{eqnarray*} \langle M_{e_1}Ue_{\gamma}, e_{\xi} \rangle &=& \langle e_{5\gamma +1}, e_{\xi} \rangle \\ &=& \widehat{\mu}(1+5\gamma - \xi) \\ &=& \langle S_1^*US_1e_{\gamma}, e_{\xi} \rangle \\ &\simeq & U|_{W_0}(\gamma, \xi) \quad \textrm{by Equation \eqref{eqn:plusone}}. \end{eqnarray*}

The result is the following commutative diagram.
\begin{center}
$\begin{CD}
L^2(\mu_{\frac14}) @>S_1>> W_0 \\
@VV M_{e_1}U V  @VV U|_{W_0}V\\
L^2(\mu_{\frac14}) @>S_1>> W_0
\end{CD}$
\end{center}
The result then follows from Theorem \ref{thm:block}.
\end{proof}

The exponentials $E(\Gamma_0)$ form an ONB for the subspace $W_0=S_1(L^2(\mu_{\frac14}))$.  We make use of the Cuntz relations on the operators $S_0$ and $S_1$ to divide $W_0$ into a direct sum of subspaces.  Let \begin{equation}\label{eqn:Vk}
V_k = S_1S_0^kS_1(L^2(\mu_{\frac14})), \qquad k=0, 1, 2, \ldots.
\end{equation}  
Each $V_k$ has ONB $E(\widetilde{\Gamma}_k)$, where we define \[ \widetilde{\Gamma}_k= 1+4^{k+1}(1+4\Gamma).\]  Using the correspondence from Equation \eqref{Eqn:SeqCor} between elements of $\Gamma$ and infinite sequences, an element $\tilde{\gamma}$ from $\widetilde{\Gamma}_k$ corresponds to a sequence where $a_0= a_{k+1} = 1$ and $a_1 = \cdots = a_k = 0$:
\[ \tilde{\gamma} \leftrightarrow (1, 0, 0, \ldots, 0, 1, a_{k+2}, a_{k+3}, \ldots, a_m, 0, 0, 0 \ldots) .\]
We see that $\widetilde{\Gamma}_k \subset \Gamma_0$ for all $k \in \mathbb{N}_0$ and that $\widetilde{\Gamma}_j$, $\widetilde{\Gamma}_k$ are pairwise mutually disjoint.  Moreover, the $\widetilde{\Gamma}_k$ sets form a partition of $\Gamma_0 \setminus \{1\}$:  
\begin{equation} \Gamma_0 = \{1\} \cup \bigsqcup_{k=0}^{\infty} \widetilde{\Gamma}_k. \end{equation}
Making use of Equation \eqref{eqn:Vk}, we have the following decomposition of $W_0$:
\begin{equation} 
W_0 = S_1(L^2(\mu)) = sp\{e_1\} \oplus \bigoplus_{k=0}^{\infty} V_k. \end{equation}

We now compute the entries of the matrix for $U|_{W_0}$.  We first observe that the majority of the blocks of the matrix are zero, and then we focus on the four row and column blocks with possibly nonzero entries.

\begin{lemma}
Let $i, j \geq 1$.  Let $\xi'\in\widetilde{\Gamma}_i = 1 + 4^{i+1}(1+4\Gamma)$ and $\gamma'\in \widetilde{\Gamma}_j=1 + 4^{j+1}(1+4\Gamma)$.  Then
\[ U|_{W_0}(\xi', \gamma') = 0,\]
so the $(\widetilde{\Gamma_i}, \widetilde{\Gamma}_j)$ block in the matrix of $U|_{W_0}$ is the zero matrix.
\end{lemma}
\begin{proof}
Let $\xi' = 1+4^{i+1}(1+4\xi)$ and $\gamma' = 1+4^{j+1}(1+4\gamma)$ for $\xi, \gamma \in \Gamma$.  Since $U|_{W_0}(\xi', \gamma') = \widehat{\mu}_{\frac14}(5\gamma' - \xi')$, we compute $5\gamma' - \xi'$:
\begin{equation*}
\begin{split}
5\gamma' - \xi'
& = 5\Bigl(1+4^{i+1}(1+4\gamma)\Bigr) - \Bigl(1+4^{j+1}(1+4\xi)\Bigr)\\
& = 5 + 5\cdot 4^{i+1}(1+4\gamma) - 1 - 4^{j+1}(1+4\xi)\\
& = 4\Bigl(1 + 4^{i}(1+4\gamma) - 4^{j}(1+4\xi) \Bigr)\\
& = 4 \cdot \textrm{odd integer} \in \mathcal{Z}(\widehat{\mu}_{\frac14}).
\end{split}
\end{equation*}   
\end{proof}

\begin{lemma}\label{lem:Umatrix} Let $\gamma'\in \Gamma_0 = 1+4\Gamma$.  The following hold for the matrix entries in the first row  of $U|_{W_0}$. 
\begin{enumerate}[\rm(a)] 
\item \rm Block ($1,1$): \:\:\it $U|_{W_0}(1,1) = 0$;
\item \rm Block ($1,\widetilde{\Gamma}_0$): \:\it If $\gamma' \in \widetilde{\Gamma}_0 = 1+4(1+4\Gamma)$, then $U|_{W_0}(1,\gamma') \neq 0$;
\item \rm Block ($1,\widetilde{\Gamma}_k$), $k > 0$: \:\it If $\gamma' \in \widetilde{\Gamma}_k = 1 + 4^{k+1}(1+4\Gamma)$ for $k>0$, then $U|_{W_0}(1,\gamma') = 0$.
\end{enumerate}
\end{lemma}
\begin{proof} 

We work with the equation $U(1,\gamma') = \widehat{\mu}_{\frac14}(5\gamma'-1)$ for $\gamma'\in \Gamma_0$.
\begin{enumerate}[(a)]
\item We have $5(1)-1 = 4 \in \mathcal{Z}(\widehat{\mu}_{\frac14})$ by \eqref{Eqn:FourierZeroes}, so $U|_{W_0}(1,1) = 0$.

\item Since $\gamma' \in \widetilde{\Gamma}_0$, we have $\gamma' = 1+4(1+4\gamma)$ for some $\gamma \in \Gamma$.  Then
\begin{equation*}
\begin{split}
5\gamma'-1 
& = 5(1+4(1+4\gamma)) - 1\\
& = 4+20+80\gamma = 4(6+20\gamma) \notin \mathcal{Z}(\widehat{\mu}_{\frac14}).
\end{split}
\end{equation*}

\item  We repeat the above computation, but now $k \geq 1$:
\begin{equation*}
\begin{split}
5\gamma'-1
& = 5(1+4^{k+1}(1+4\gamma)) - 1\\
&=  4+5\cdot 4^{k+1}+ 5 \cdot 4^{k+2}\gamma \\ 
&= 4( 1+5 \cdot 4^{k}+ 5 \cdot 4^{k+1}\gamma) \\ 
&= 4 \cdot \textrm{odd integer} \in \mathcal{Z}(\widehat{\mu}_{\frac14}).
\end{split}
\end{equation*}  This proves  $\widehat{\mu}_{\frac14}(5\gamma'-1) = 0$.
\end{enumerate}
\end{proof}

We now turn to the first column of the matrix of $U|_{W_0}$.
\begin{lemma}
The following hold for the matrix entries in the first column of $U|_{W_0}$. 
\begin{enumerate}[\rm(a)] 
\item \rm Block ($\widetilde{\Gamma}_0, 1$): \:\it This block contains nonzero entries.
\item \rm Block ($\widetilde{\Gamma}_k, 1$), $k > 0$: \:\it Each of the blocks $(\widetilde{\Gamma}_k, 1)$ is zero.
\end{enumerate}
\end{lemma}
\begin{proof}
Let $\xi' = 1+4^{k+1}(1+4\xi)$ , where $\xi\in\Gamma$.
\begin{enumerate}[(a)] 
\item Suppose $k = 0$.  We compute $5 -\xi'$:
\begin{equation*}
\begin{split}
5 - (1+4(1+4\xi))
& = 4 - 4(1+4\xi)\\
& = 4 (1 - (1+4\xi))\\
& = -4^2\xi. 
\end{split}
\end{equation*}
If $\xi\in\Gamma$ is even, which it can be ($\xi = 4, 16, 20, \ldots$), the associated matrix entry is nonzero.

\item We compute $5 - \xi' = 5 - (1+4^{i+1}(1+4\xi))$, where $k \geq 1$:
\begin{equation*}
\begin{split}
5 - (1+4^{k+1}(1+4\xi))
& = 4 - 4^{k+1}(1+4\xi)\\
& = 4 (1 - 4^k(1+4\xi))\\
& = 4 \cdot \textrm{odd integer} \in \mathcal{Z}(\widehat{\mu}_{\frac14}).
\end{split}
\end{equation*}
\end{enumerate}
\end{proof}
Finally, we examine the $(\widetilde{\Gamma}_k, \widetilde{\Gamma}_0)$ and $(\widetilde{\Gamma}_0, \widetilde{\Gamma}_k)$ blocks in the matrix for $U|_{W_0}$.
\begin{lemma}\label{lem:UW0matrix} Let $k \geq 1$. 
\begin{enumerate}[\rm(a)] 
\item \rm Block ($ \widetilde{\Gamma}_0, \widetilde{\Gamma}_0$): \:\:\it This block is zero.
\item \rm Block ($ \widetilde{\Gamma}_0, \widetilde{\Gamma}_k$): \:\it  This block contains nonzero entries.
\item \rm Block ($\widetilde{\Gamma}_k, \widetilde{\Gamma}_0$): \:\it  This block contains nonzero entries.
\end{enumerate}
\end{lemma}
\begin{proof} Again, we use the formula $U|_{W_0}(\xi', \gamma') = \widehat{\mu}_{\frac14}(5\gamma'-\xi')$. \begin{enumerate}[(a)]
\item  Let $\xi' = 1+4(1+4\xi)$ and $\gamma' = 1+4(1+4\gamma)$, where $\xi, \gamma\in\Gamma$.  Then
\begin{equation*}
\begin{split}
5\gamma' - \xi'
& = 5\Bigl(1+4(1+4\gamma) \Bigr) - 1 - 4(1+4\xi)\\
& = 5 + 5\cdot 4(1+4\gamma)- 1 - 4(1+4\xi)\\
& = 4\Bigl(1 + 5(1+4\gamma) - 1-4\xi)\Bigr)\\
& = 4\Bigl(5+20\gamma-4\xi\Bigr)\in \mathcal{Z}(\mu_{\frac14}).\\
\end{split}
\end{equation*}
\item  Let $\xi' = 1+4(1+4\xi)$ and $\gamma' = 1+4^{k+1}(1+4\gamma)$, where $\xi, \gamma\in\Gamma$ and $k\geq 1$.  Then
\begin{equation*}
\begin{split}
5\gamma' - \xi'
& = 5\Bigl(1+4^{k+1}(1+4\gamma)\Bigr)- 1 - 4(1+4\xi)\\
& = 5 + 5\cdot 4^{k+1}(1+4\gamma)- 1 - 4(1+4\xi)\\
&= 4\Bigl(1 + 4^{k}(1+4\gamma) - 1 - 4\xi \Bigr)\\
& = 4\Bigl(4^{k}(1+4\gamma) - 4\xi \Bigr)\\
& = 4^2\Bigl(4^{k-1}(1+4\gamma) - \xi \Bigr)\\
\end{split}
\end{equation*}
which is not necessarily a zero of $\widehat{\mu}_{\frac14}$.  For example, if $k$ and $\gamma$ are fixed, then $4^{k-1}(1+4\gamma)\in \Gamma$.  Set $\xi = 4^{k-1}(1+4\gamma)$.  Then $5\gamma' - \xi' = 0$, which does not belong to the zero set of $\widehat{\mu}_{\frac14}$.
\item Let $\xi' = 1+4^{k+1}(1+4\xi)$ and $\gamma' = 1+4(1+4\gamma)$, where $\xi, \gamma\in\Gamma$ and $k\geq 1$.  Then
\begin{equation*}
\begin{split}
5\gamma' - \xi'
& = 5\Bigl(1+4(1+4\gamma)\Bigr)- 1 - 4^{k+1}(1+4\xi)\\
& = 5 + 5\cdot 4(1+4\gamma)- 1 - 4^{k+1}(1+4\xi)\\
& = 4\Bigl(1 + 5(1+4\gamma)- 4^{k}(1+4\xi)  \Bigr)\\
& =  4\Bigl(6+ 20\gamma- 4^{k}(1+4\xi)  \Bigr),\\
\end{split}
\end{equation*}
which cannot belong to the zero set of $\widehat{\mu}_{\frac14}$; $\frac{1}{4}(5\gamma' - \xi')$ is congruent to $2 \pmod{4}$.
\end{enumerate} 
\end{proof}

The results from our previous lemmas combine to give a block structure to the matrix for the operator $U|_{W_0}$.

\begin{theorem}\label{Thm:SpecialCase}   The matrix for $U|_{W_0}$ has the following block form, where the index set $\Gamma_0$ for the ONB of $W_0$ is ordered by $\{1\}\cup \bigsqcup \widetilde{\Gamma}_k$.  The notation $\ast$ indicates an entry or block  that is not necessarily zero and $0$ indicates an entry or block that contains all zeros.

\begin{equation}\label{Eqn:MatrixU_W_0}
U|_{W_0} \simeq  \begin{tabular}{c||c|c|c|c|c|c}
                         & \:\:$1$\:\: & $\widetilde{\Gamma}_0$  & $\widetilde{\Gamma}_1$ & $\widetilde{\Gamma}_2$ &  $\widetilde{\Gamma}_3$ & $\cdots$\\
\hline
\hline
$1$ & $0$ &$\ast$ & $0$ & $0$ & $0$ &$\cdots$\\ 
\hline $\widetilde{\Gamma}_0$  & $\ast$ & $0$ &    $\ast$        &     $\ast$   &  $\ast$   &       $\cdots$ \\
\hline 
$\widetilde{\Gamma}_1$  & $0$ &  $\ast$          &$0$ &       $0$   &          $0$   &      $\cdots$ \\
\hline
$\widetilde{\Gamma}_2$  & $0$ & $\ast$           &   $0$          &        $0$   &         $0$   &       $\cdots$ \\

\hline
$\widetilde{\Gamma}_3$  & $0$ & $\ast$           &   $0$          &        $0$   &         $0$   &       $\cdots$ \\
\hline

                        & \vdots &  $\vdots$   &    $\vdots$   &      $\vdots$                &  $\vdots$                &  $\ddots$ \\ 
\end{tabular}  \end{equation}

\end{theorem}

\bibliographystyle{alpha}
\bibliography{U_structure}

\begin{thebibliography}{DHSW11}

\bibitem[BJ99]{BrJo99}
Ola Bratteli and Palle E.~T. Jorgensen.
\newblock Iterated function systems and permutation representations of the
  {C}untz algebra.
\newblock {\em Mem. Amer. Math. Soc.}, 139(663):x+89, 1999.

\bibitem[Cun77]{Cun77}
Joachim Cuntz.
\newblock Simple {$C\sp*$}-algebras generated by isometries.
\newblock {\em Comm. Math. Phys.}, 57(2):173--185, 1977.

\bibitem[DHJ09]{DHJ09}
Dorin~Ervin Dutkay, Deguang Han, and Palle E.~T. Jorgensen.
\newblock Orthogonal exponentials, translations, and {B}ohr completions.
\newblock {\em J. Funct. Anal.}, 257(9):2999--3019, 2009.

\bibitem[DHS09]{DHS09}
Dorin~Ervin Dutkay, Deguang Han, and Qiyu Sun.
\newblock On the spectra of a {C}antor measure.
\newblock {\em Adv. Math.}, 221(1):251--276, 2009.

\bibitem[DHSW11]{DHSW11}
Dorin~Ervin Dutkay, Deguang Han, Qiyu Sun, and Eric Weber.
\newblock On the {B}eurling dimension of exponential frames.
\newblock {\em Adv. Math.}, 226(1):285--297, 2011.

\bibitem[DJ09]{DJ09}
Dorin~Ervin Dutkay and Palle E.~T. Jorgensen.
\newblock Fourier duality for fractal measures with affine scales.
\newblock ar{X}iv:0911.1070, 2009.

\bibitem[DJ11]{DJ10}
Dorin~Ervin Dutkay and Palle E.~T. Jorgensen.
\newblock Spectral measures and {C}untz algebras.
\newblock ar{X}iv:1001.4565v1, 2011.

\bibitem[Erd39]{Erd39}
Paul Erd{\"o}s.
\newblock On a family of symmetric {B}ernoulli convolutions.
\newblock {\em Amer. J. Math.}, 61:974--976, 1939.

\bibitem[Fug74]{Fug74}
Bent Fuglede.
\newblock Commuting self-adjoint partial differential operators and a group
  theoretic problem.
\newblock {\em J. Functional Analysis}, 16:101--121, 1974.

\bibitem[GN98]{GaNa98}
Jean-Pierre Gabardo and M.~Zuhair Nashed.
\newblock Nonuniform multiresolution analyses and spectral pairs.
\newblock {\em J. Funct. Anal.}, 158(1):209--241, 1998.

\bibitem[GY06]{GaYu06}
Jean-Pierre Gabardo and Xiaojiang Yu.
\newblock Wavelets associated with nonuniform multiresolution analyses and
  one-dimensional spectral pairs.
\newblock {\em J. Math. Anal. Appl.}, 323(2):798--817, 2006.

\bibitem[HL08]{HuLa08}
Tian-You Hu and Ka-Sing Lau.
\newblock Spectral property of the {B}ernoulli convolutions.
\newblock {\em Adv. Math.}, 219(2):554--567, 2008.

\bibitem[Hut81]{Hut81}
John~E. Hutchinson.
\newblock Fractals and self-similarity.
\newblock {\em Indiana Univ. Math. J.}, 30(5):713--747, 1981.

\bibitem[JKS07]{JKS07b}
Palle E.~T. Jorgensen, Keri~A. Kornelson, and Karen~L. Shuman.
\newblock Harmonic analysis of iterated function systems with overlap.
\newblock {\em J. Math. Phys.}, 48(8):083511, 35, 2007.

\bibitem[JKS08]{JKS08}
Palle E.~T. Jorgensen, Keri~A. Kornelson, and Karen~L. Shuman.
\newblock Orthogonal exponentials for {B}ernoulli iterated function systems.
\newblock In {\em Representations, Wavelets, and Frames: A celebration of the
  mathematical work of Lawrence W. Baggett}, Appl. Num. Harm. Anal., pages
  217--237. Birkh\"{a}user, Boston, MA, 2008.

\bibitem[JKS11]{JKS11a}
Palle E.~T. Jorgensen, Keri~A. Kornelson, and Karen~L. Shuman.
\newblock Families of spectral sets for {B}ernoulli convolutions.
\newblock {\em J. Fourier Anal. Appl.}, 17(3):431--456, 2011.

\bibitem[JP98]{JoPe98}
Palle E.~T. Jorgensen and Steen Pedersen.
\newblock Dense analytic subspaces in fractal {$L\sp 2$}-spaces.
\newblock {\em J. Anal. Math.}, 75:185--228, 1998.

\bibitem[JR95]{JR95}
Aimee Johnson and Daniel~J. Rudolph.
\newblock Convergence under {$\times_q$} of {$\times_p$} invariant measures on
  the circle.
\newblock {\em Adv. Math.}, 115(1):117--140, 1995.

\bibitem[Li07]{Li07}
Jian-Lin Li.
\newblock Spectral self-affine measures in {$\Bbb R\sp N$}.
\newblock {\em Proc. Edinb. Math. Soc. (2)}, 50(1):197--215, 2007.

\bibitem[{\L}W02]{LaWa02}
Izabella {\L}aba and Yang Wang.
\newblock On spectral {C}antor measures.
\newblock {\em J. Funct. Anal.}, 193(2):409--420, 2002.

\bibitem[PS96]{PeSo96}
Yuval Peres and Boris Solomyak.
\newblock Absolute continuity of {B}ernoulli convolutions, a simple proof.
\newblock {\em Math. Res. Lett.}, 3(2):231--239, 1996.

\bibitem[PSS00]{PSS00}
Yuval Peres, Wilhelm Schlag, and Boris Solomyak.
\newblock Sixty years of {B}ernoulli convolutions.
\newblock In {\em Fractal geometry and stochastics, {II}
  ({G}reifswald/{K}oserow, 1998)}, volume~46 of {\em Progr. Probab.}, pages
  39--65. Birkh\"auser, Basel, 2000.

\bibitem[PW01]{PeWa01}
Steen Pedersen and Yang Wang.
\newblock Universal spectra, universal tiling sets and the spectral set
  conjecture.
\newblock {\em Math. Scand.}, 88(2):246--256, 2001.

\bibitem[Sid03]{Sid03}
N.~Sidorov.
\newblock Ergodic-theoretic properties of certain {B}ernoulli convolutions.
\newblock {\em Acta Math. Hungar.}, 101(4):345--355, 2003.

\bibitem[Tao04]{Tao04}
Terence Tao.
\newblock Fuglede's conjecture is false in 5 and higher dimensions.
\newblock {\em Math. Res. Lett.}, 11(2-3):251--258, 2004.

\bibitem[Wan02]{Wan02}
Yang Wang.
\newblock Wavelets, tiling, and spectral sets.
\newblock {\em Duke Math. J.}, 114(1):43--57, 2002.

\end{thebibliography}

\end{document}